\documentclass[11pt]{article}
\makeatletter

\newcommand{\Rmnum}[1]{\expandafter\@slowromancap\romannumeral #1@}
\makeatother
\usepackage{}
\usepackage{mathrsfs}
\usepackage{amsfonts}
\usepackage{amsfonts,amssymb,mathrsfs}
\usepackage{amsmath,amscd}
\usepackage[dvips]{graphicx}
\usepackage[all]{xy}
\usepackage{graphicx}
\usepackage{fancyhdr}
\usepackage{mathptm,pslatex}
\usepackage{amsthm}

\begin{document}

\sloppy
\newtheorem{Def}{Definition}[section]
\newtheorem{Bsp}{Example}[section]
\newtheorem{Prop}[Def]{Proposition}
\newtheorem{Theo}[Def]{Theorem}
\newtheorem{Lem}[Def]{Lemma}
\newtheorem{Koro}[Def]{Corollary}
\theoremstyle{definition}
\newtheorem{Rem}[Def]{Remark}

\newcommand{\add}{{\rm add}}
\newcommand{\gd}{{\rm gl.dim }}
\newcommand{\dm}{{\rm dom.dim }}
\newcommand{\E}{{\rm E}}
\newcommand{\Mor}{{\rm Morph}}
\newcommand{\End}{{\rm End}}
\newcommand{\ind}{{\rm ind}}
\newcommand{\rsd}{{\rm res.dim}}
\newcommand{\rd} {{\rm rep.dim}}
\newcommand{\ol}{\overline}
\newcommand{\rad}{{\rm rad}}
\newcommand{\soc}{{\rm soc}}
\renewcommand{\top}{{\rm top}}
\newcommand{\pd}{{\rm proj.dim}}
\newcommand{\rpd}{{\rm re.proj.dim}}
\newcommand{\id}{{\rm inj.dim}}
\newcommand{\Fac}{{\rm Fac}}
\newcommand{\fd} {{\rm fin.dim }}
\newcommand{\rfd} {{\rm re.fin.dim }}
\newcommand{\DTr}{{\rm DTr}}
\newcommand{\cpx}[1]{#1^{\bullet}}
\newcommand{\D}[1]{{\mathscr D}(#1)}
\newcommand{\Dz}[1]{{\mathscr D}^+(#1)}
\newcommand{\Df}[1]{{\mathscr D}^-(#1)}
\newcommand{\Db}[1]{{\mathscr D}^b(#1)}
\newcommand{\C}[1]{{\mathscr C}(#1)}
\newcommand{\Cz}[1]{{\mathscr C}^+(#1)}
\newcommand{\Cf}[1]{{\mathscr C}^-(#1)}
\newcommand{\Cb}[1]{{\mathscr C}^b(#1)}
\newcommand{\K}[1]{{\mathscr K}(#1)}
\newcommand{\Kz}[1]{{\mathscr K}^+(#1)}
\newcommand{\Kf}[1]{{\mathscr  K}^-(#1)}
\newcommand{\Kb}[1]{{\mathscr K}^b(#1)}
%\stackrel{\sim}
\newcommand{\Modcat}{\ensuremath{\mbox{{\rm -Mod}}}}
\newcommand{\modcat}{\ensuremath{\mbox{{\rm -mod}}}}
\newcommand{\stmodcat}[1]{#1\mbox{{\rm -{\underline{mod}}}}}
\newcommand{\pmodcat}[1]{#1\mbox{{\rm -proj}}}
\newcommand{\imodcat}[1]{#1\mbox{{\rm -inj}}}
\newcommand{\opp}{^{\rm op}}
\newcommand{\otimesL}{\otimes^{\rm\bf L}}
\newcommand{\rHom}{{\rm\bf R}{\rm Hom}}
\newcommand{\projdim}{\pd}
\newcommand{\Hom}{{\rm Hom}}
\newcommand{\Coker}{{\rm coker}\,\,}
\newcommand{ \Ker  }{{\rm Ker}\,\,}
\newcommand{ \Img  }{{\rm Im}\,\,}
\newcommand{\Ext}{{\rm Ext}}
\newcommand{\Tor}{{\rm Tor}}
\newcommand{\StHom}{{\rm \underline{Hom} \, }}

\newcommand{\gm}{{\rm _{\Gamma_M}}}
\newcommand{\gmr}{{\rm _{\Gamma_M^R}}}

\def\vez{\varepsilon}\def\bz{\bigoplus}  \def\sz {\oplus}
\def\epa{\xrightarrow} \def\inja{\hookrightarrow}

\newcommand{\lra}{\longrightarrow}
\newcommand{\lraf}[1]{\stackrel{#1}{\lra}}
\newcommand{\ra}{\rightarrow}
\newcommand{\dk}{{\rm dim_{_{k}}}}

{\Large \bf
\begin{center}
Finitistic dimension conjecture and extensions of algebras
\end{center}}
\medskip

\centerline{\bf {Shufeng Guo$^{a, b}$}}

\begin{center} $^{a}$ Faculty of Science, Guilin University of Aerospace Technology, 541004 Guilin, \\People's Republic of China
\end{center}

\begin{center} $^{b}$ School of Mathematical Sciences, Capital Normal University, 100048 Beijing, \\People's Republic of  China
\end{center}

\renewcommand{\thefootnote}{\alph{footnote}}
\setcounter{footnote}{-1} \footnote{E-mail address: guoshufeng132@126.com}
%\date{}

\begin{abstract}
An extension of algebras is a homomorphism of algebras preserving identities. We use extensions of algebras to study the finitistic dimension conjecture over Artin algebras. Let $f: B \to A$ be an extension of Artin algebras. We denote by $\fd(f)$ the relative finitistic dimension of $f$, which is defined to be the supremum of relative projective dimensions of finitely generated left $A$-modules of finite projective dimension. We prove that, if $B$ is representation-finite and $\fd(f)\leq 1$, then $A$ has finite finitistic dimension. For the case of $\fd(f)> 1$, we give a sufficient condition for $A$ with finite finitistic dimension. Also, we prove the following result: Let $I$, $J$, $K$ be three ideals of an Artin algebra $A$ such that $IJK=0$ and $K\supseteq \rad(A)$. If both $A/I$ and $A/J$ are $A$-syzygy-finite, then the finitistic dimension of $A$ is finite.

\end{abstract}

\medskip


\section{Introduction}
Let $A$ be an Artin algebra. The finitistic dimension of $A$ is defined to be the supremum of projective dimensions of finitely generated left $A$-modules having finite projective dimension. The famous finitistic dimension conjecture says that the finitistic dimension of any Artin algebra is finite (see \cite[Conjecture 11, pp. 410]{ARS} or \cite{B}). It is 57 years old and remains open to date. It is worth noting that the finitistic dimension conjecture is very closely related to many homological conjectures in the representation theory of algebras, such as strong Nakayama conjecture, generalized Nakayama conjecture, Nakayama conjecture, Wakamatsu tilting conjecture and Gorenstein symmetry conjecture. If the finitistic dimension conjecture holds, then so do the above conjectures (\cite{ARE, Y}). However, there are a few cases for which this conjecture is verified to be true (see, for example, \cite{GZ, GKK, EHIS}). In general, this conjecture seems to be far from being solved.

Recently, the work of Xi in \cite{X1, X2} shows that the finitistic dimension conjecture can be reduced to comparing finitistic dimensions of two algebras in an extension. The basic idea is as follows: let $B$ and $A$ be Artin algebras, and $f: B \to A$ a homomorphism of algebras satisfying some certain conditions. If one of them has finite finitistic dimension, is the finitistic dimension of the other finite? From on the other hand of view, it is reasonable to study the finitistic dimension conjecture by extensions of algebras. In fact, we have known that some classes of algebras have finite finitistic dimension, so we use them to obtain more classes of algebras with finite finitistic dimension by means of extension. In literatures, we have already seen some interesting results concerning this direction (see \cite{EHIS, X1, X2, XX, Wjq, Wjq2, WX}). In this note, we shall continue to study the above question.

Different from the usual consideration (see, for example, \cite{X1, X2, XX}), where one often uses the information on $A$ to get the information on $B$, we use some relative homological dimension to control the extension $f: B \to A$ and employ the finitistic dimension of $B$ to study that of $A$. Here, the relative finitistic dimension of $f$, denoted by $\fd(f)$, is defined to be the supremum of relative projective dimensions of finitely generated left $A$-modules of finite projective dimension. We get the following result, which generalizes the result of E. L. Green in \cite[Theorem 1.5]{G}.

\begin{Theo} Let $B$ and $A$ be Artin algebras with $B$ representation-finite. Suppose that $\varphi: B \to A$  is a homomorphism of algebras preserving identities. Then:

$(1)$ If $\fd(\varphi)\leq 1$, then $A$ has finite finitistic dimension.

$(2)$ If $2\leq\fd(\varphi)<\infty$ and if, for any $A$-module $X$ with finite projective dimension, $_{A}A\otimes_{B}X$ has finite projective dimension, then $A$ has finite finitistic dimension. \label{thm1}
\end{Theo}

\medskip
In Theorem \ref{thm1}, we use the finitistic dimension of $B$ to describe that of $A$. In the following, for an extension $f: B \to A$, we shall employ the finiteness of the finitistic dimension of $A$ to approach that of $B$. On the one hand, we establish the relationship between the finiteness of finitistic dimensions of quotient algebras and given algebras, and obtain the following result, which recovers many known results in literatures, for example, \cite[Theorem 3.2, Lemma 3.6, Corollary 3.8]{X1}, \cite[Theorem 3.1, Corollary 3.2, Corollary 3.3, Proposition 3.5]{Wjq2}, the result in \cite{W} and so on. For unexplained notions in the following result, we refer to Section \ref{pre}.

\begin{Theo} Let $A$ be an Artin algebra and let $I$, $J$, $K$ be three ideals of $A$ such that $IJK=0$ and $K\supseteq \rad(A)$. If both $A/I$ and $A/J$ are $A$-syzygy-finite, then the finitistic dimension of $A$ is finite.  \label{thm2}
\end{Theo}

\medskip
On the other hand, we consider left idealized extensions to study the finitistic dimension conjecture, and get the following.

\begin{Prop} Let

$$B=A_{0}\subseteq A_{1}\subseteq\cdots  \subseteq A_{s-1}\subseteq A_{s}=A$$

\noindent be a chain of subalgebras of an Artin algebra $A$ such that $\rad (A_{i-1})$ is a left ideal of $A_{i}$ for all $1\leq i\leq s$ with $s$ being a positive integer and that $A$ is 1-syzygy-finite. Then $\fd(B)<\infty$ provided one of the following conditions is satisfied.

$(1)$ $B/\rad (A_{s-1})\cdots \rad (A_{1})\rad (A_{0})$ is $B$-syzygy-finite (for example, $B/\rad (A_{s-1})\cdots \rad (A_{1})\rad(A_{0})$ is representation-finite).

$(2)$ $A_{1}/\rad (A_{s-1})\cdots \rad (A_{1})$ is $B$-syzygy-finite (for example, $A_{1}/\rad (A_{s-1})\cdots \rad (A_{1})$ is representation-finite). \label{propo1}
\end{Prop}

\medskip
Remark that Proposition \ref{propo1} recovers \cite[Theorem 3.1]{X1} if we take $s=1$, and reobtain \cite[Theorem 4.5]{X1} if we take $s=2$.

The paper is organized as follows. In Section \ref{pre} we recall some definitions and basic results which are need in the paper.
We prove Theorem \ref{thm1} in Section \ref{rfd} and give a proof of Theorem \ref{thm2} in Section \ref{qfd}. In the last section we use left idealized extensions to study the finitistic dimension conjecture and prove Proposition \ref{propo1}.

\section{Preliminaries\label{pre}}
In this section, we shall fix some notations, and recall some
definitions and basic results which are needed in the proofs of our main results. Throughout this paper, unless stated otherwise, all the algebras considered are Artin $R$-algebras, where $R$ is assumed to be a commutative Artin ring, and all the modules considered are finitely generated left modules over Artin algebras, so that all the homological dimensions will be assumed to be in the category of finitely generated modules.

Let $A$ be an Artin algebra. We denote by $A\modcat$ the category of all finitely
generated left $A$-modules, and by $\rad(A)$ the Jacobson radical of $A$. Given an $A$-module $M$,
we denote by $\pd(_{A}M)$ the projective dimension of $M$, by $\Omega_{A}^{i}(M)$ the $i$-th
syzygy of $M$ (we set $\Omega_{A}^{0}(M):=M$), and by $\add(_{A}M)$ the full subcategory of $A\modcat$ consisting of all direct summands
of finite direct sums of copies of $M$.

Now let us recall some definitions concerning Artin algebras. $A$ is called representation-finite if there is only finitely many nonisomorphic indecomposable $A$-modules in $A\modcat$. The finitistic dimension of $A$, denoted by
$\fd(A)$, is defined as
$$\begin{array}{rl} \fd(A)=\mbox{Sup} \{\pd(_AM) \mid M\in A\modcat  \mbox{ and }  \pd(_AM)<\infty \}. \end{array}$$

\noindent And the global dimension of $A$, denoted by $\gd(A)$, is defined as
$$\begin{array}{rl} \gd(A)=\mbox{Sup} \{\pd(_AM) \mid M\in A\modcat  \}.  \end{array}$$

Let $\mathcal{C}$ be a subcategory of $A\modcat$ and $m$ a natural number. We set
$$ \begin{array}{rl} \Omega_{A}^{m}(\mathcal{C}):=\{\Omega_{A}^{m}(X) \mid X\in \mathcal{C}  \}. \end{array}$$

\noindent $\mathcal{C}$ is said to be $m$-$A$-syzygy-finite, or simply $m$-syzygy-finite if there is no confusion, if
the number of non-isomorphic indecomposable direct summands of objects in $\Omega_{A}^{m}(\mathcal{C})$ is finite, that is, there is an $A$-module $N$
such that $\Omega_{A}^{m}(\mathcal{C})\subseteq \add(_{A}N)$. Furthermore,
we say that $\mathcal{C}$ is $(A)$-syzygy-finite if there is some natural number $n$ such that $\mathcal{C}$ is $n$-$(A)$-syzygy-finite.
If $A\modcat$ is syzygy-finite, then we also say that $A$ is syzygy-finite. Let $C$ be a second Artin algebra and $f: A\to C$ a homomorphism
of algebras preserving identities. Clearly, every $C$-module can be regarded as an $A$-module in the natural way, and every $C$-homomorphism
can be viewed as an $A$-homomorphism. This means that $C\modcat$ is a subcategory of $A\modcat$. If $C\modcat$ is $A$-syzygy-finite, then we also say that
$C$ is $A$-syzygy-finite. Note that if $C$ is representation-finite, then $C$ is $A$-syzygy-finite.

Next we give the definition and basic properties of Igusa-Todorov function. We denote by $K_{0}(A)$ the free abelian group generated by the isomorphism classes $[M]$ of modules $M$ in $A\modcat$. Let $K(A)$ be the factor group of $K_{0}(A)$ modulo the following relations:

\medskip
(1) $[Y]=[X]+[Z]$ if $Y\simeq X\oplus Z$;

(2) $[P]=0$ if $P$ is projective.

\medskip
\noindent Then $K(A)$ is also the free abelian group with basis the isomorphism classes of indecomposable non-projective $A$-modules in $A\modcat$. Igusa and Todorov in \cite{IT} introduced a function $\Psi: A\modcat \rightarrow \mathbb{N}$ on this abelian group, which is defined on the objects of $A\modcat$ and takes values of non-negative integers. We call it the Igusa-Todorov function. It follows from \cite{IT} that, for any Artin algebra $A$, the Igusa-Todorov function always exists. For the convenience of the reader, we give the basic properties of Igusa-Todorov function as follows.

\begin{Lem} (\cite{IT}) Let $A$ be an Artin algebra and $\Psi$ be the {\rm Igusa-Todorov} function. Then the following are true.

(1) For any $A$-module $M$, if $M$ has finite projective dimension, then $\Psi(M)=\pd(_{A}M)$.

(2) If $0\ra X\ra Y \ra Z\ra 0$ is an exact sequence in $A\modcat$ with $\pd(Z)< \infty$, then $\pd(Z)\leq\Psi(X\oplus Y)+1$.

(3) If $0\ra X\ra Y \ra Z\ra 0$ is an exact sequence in $A\modcat$ with $\pd(Y)< \infty$, then $\pd(Y)\leq\Psi(\Omega(X)\oplus\Omega^{2}(Z))+2$.

(4) If $0\ra X\ra Y \ra Z\ra 0$ is an exact sequence in $A\modcat$ with $\pd(X)< \infty$, then $\pd(X)\leq\Psi(\Omega(Y\oplus Z))+1$.
\label{lem2.1}
\end{Lem}

\medskip
Finally, we shall recall some definitions and basic facts on relative homological algebra.
Let $B$ and $A$ be Artin algebras, and $f: B\to A$ a homomorphism of algebras preserving identities.
Then we say that $f$ is an extension. Clearly, every $A$-module can
be regarded as a $B$-module via $f$ in the natural way. An exact sequence in $A\modcat$

$$\cdots \lra M_{i+1}\lra M_{i} \lraf {t_{i}} M_{i-1} \lra\cdots$$
is called $(A, B)$-exact if there are $B$-homomorphisms $h_{i}: M_{i} \to M_{i+1}$ such that $t_{i}=t_{i}h_{i-1}t_{i}$ for all $i$.
It is very easily checked that the definition is equivalent to that introduced in \cite{H1956}.

Let $X$ be an $A$-module. $X$ is said to be $(A, B)$-projective, or relatively projective over $B$,
if $X$ is an $A$-direct summand of $A\otimes_{B}X$. For the equivalent conditions of relatively projective modules, we refer the reader to \cite[ pp. 202, Proposition 3.6]{ARS} and \cite{Hirata, T}. We denote by $\mathscr{P}(A, B)$ the full subcategory of $A\modcat$ consisting
of all $(A, B)$-projective $A$-modules. Note that $\mathscr{P}(A, B)$ is functorially finite in $A\modcat$ (see \cite{KP}).

Given an $A$-module $X$, an $(A, B)$-projective resolution of $_AX$ is defined to be an $(A, B)$-exact sequence

$$\cdots\lra P_n\lra P_{n-1}\lra\cdots\lra P_1\lra P_0\lra X\lra 0$$  in which $P_i\in \mathscr{P}(A, B)$ for each $i$.

The relative projective dimension of  $_AX$ , denoted by $\rpd(_AX)$, is defined as

$$\begin{array}{rl} \rpd(_AX)=\mbox{Inf} &\{\mbox{ n } \mid 0\ra P_n\ra P_{n-1}\ra\cdots\ra P_1\ra P_0\ra X\ra 0 \\
  &\qquad   \mbox{ is an } (A, B)\mbox{-projective resolution of }_AX\}.  \end{array}$$

\noindent If such an exact sequence does not exist, we say that the relative projective dimension of $_AX$ is infinite.
The relative global dimension of the extension $f$, denoted by $\gd(f)$, is defined as

 $$\begin{array}{rl} \gd(f)=\mbox{Sup} \{\rpd(_AX) \mid X\in A\modcat  \},  \end{array}$$

\noindent while the relative finitistic dimension of $f$, denoted by $\fd (f)$, is defined as
$$\begin{array}{rl} \fd(f)=\mbox{Sup} \{\rpd(_AX) \mid X\in A\modcat  \mbox{ and }  \pd(_AX)<\infty \}. \end{array}$$

\noindent Clearly, $\fd (f) \leq \gd (f)$. In particular, if $\gd(A)<\infty$, then $\fd (f) = \gd (f)$.
Xi and Xu also defined in \cite{XX} some relative finitistic dimension of $f$ to be the supremum of relative projective dimensions of finitely generated left $A$-modules with finite relative projective dimension, denoted by $\rfd (f)$.
Note that, if $\gd(B)<\infty$ and $A_{B}$ is projective, then $\rfd (f)\leq \fd (f)\leq \gd (f)$ by \cite[Theorem 1]{H1958}.

The following result is a consequence of Generalized Schanuel's Lemma in \cite{T} by induction.

\begin{Lem} \label{lemma2.2}Let $f: B\to A$ be an extension of Artin algebras.
Suppose that

$$0\lra M\lra P_{n-1}\lra P_{n-2}\lra \cdots \lra P_{0}\lra X\lra 0$$

\noindent and

$$0\lra N\lra Q_{n-1}\lra Q_{n-2}\lra \cdots \lra Q_{0}\lra X\lra 0$$

\noindent are $(A, B)$-exact sequences in which all $P_{i}$ and $Q_{i}$ are $(A, B)$-projective for $0\leq i\leq n-1$ with $n$ being a positive integer. Then we have an isomorphism

$$ M \oplus Q_{n-1}\oplus P_{n-2} \oplus \cdots \oplus C \backsimeq N\oplus P_{n-1}\oplus Q_{n-2} \oplus \cdots \oplus C' $$

\noindent as $A$-modules, where $C= P_{0}$ and $C'= Q_{0}$ if $n$ is an even number, $C= Q_{0}$ and $ C'= P_{0}$ if $n$ is an odd number. \end{Lem}

\section{Relative finitistic dimensions and finitistic dimensions \label{rfd}}

In this section, we employ the relative finitistic dimension to control an extension $\varphi: B \to A$ and use the finitistic dimension of $B$ to approach the finiteness of the finitistic dimension of $A$. Concretely, we consider the case where $B$ is of finite representation type and give a proof of Theorem \ref{thm1}.

The main result in this section is based on the following observation.

\begin{Lem} Let $\varphi: B\to A$ be an extension of Artin algebras. If $B$ is representation-finite, then so is $\mathscr{P}(A, B)$.
\label{lem3.1}
\end{Lem}

{\bf Proof.} Let $Y$ be an $A$-module in $\mathscr{P}(A, B)$. Then $Y$ is an $A$-direct summand of $A\otimes_{B}Y$. So it follows from the proof of \cite[pp. 200, Lemma 3.1]{ARS}. $\square$

\medskip
Now we can prove Theorem \ref{thm1}.

\medskip
{\bf Proof of Theorem \ref{thm1}.} Since $B$ is representation-finite, we have $\mathscr{P}(A, B)$ is also representation-finite by Lemma \ref{lem3.1}, so that we may assume that $\{Q_{1},\,Q_{2},\,\cdots,\,Q_{m}\}$ is a complete list of non-isomorphic indecomposable $(A,\,B)$-projective $A$-modules. Let $X$ be an $A$-module with finite projective dimension.

$(1)$ If $\fd(\varphi)\leq 1$, then $\rpd(_{A}X)\leq 1$, and hence $_{A}X$ has an $(A,\,B)$-projective resolution of length 1:
$$0\lra P_1\lra P_{0}\lra X\lra 0.$$

\noindent Since $_{A}P_{1}\oplus P_{0}\in \mathscr{P}(A, B)$, we can write
$_{A}P_{1}\oplus P_{0} =\oplus_{j=1}^{m}Q_{j}^{s_{j}}$, where $s_{j}$ is a non-negative integer for each $j$.
Now we bound the projective dimension of $_{A}X$:

$$\begin{array}{rl} \pd (_{A}X)&\leq \Psi(P_{1}\oplus P_{0})+1\\
&\\
&=\Psi(\oplus_{j=1}^{m}Q_{j}^{s_{j}})+1\\
&\\
&\leq \Psi(\oplus_{j=1}^{m}Q_{j})+1,

\end{array}$$

\medskip
\noindent where $\Psi$ is the {\rm Igusa-Todorov} function. Thus $\fd(A)$ is upper bounded by $\Psi(\oplus_{j=1}^{m}Q_{j})+1$ and $\fd(A)<\infty$.

$(2)$ If $2\leq\fd(\varphi)=n<\infty$, then by definition $\rpd(_{A}X)\leq n$, so $_{A}X$ has an $(A,\,B)$-projective resolution of length $n$. Consider the standard relative projective resolution of $_{A}X$

$$\cdots\lra C_n\lraf{\delta_{n}} C_{n-1}\lra\cdots\lra C_1\lra C_0\lra X\lra 0,$$

\noindent where $ C_0= A\otimes_{B}X$ and $C_{i}= A\otimes_{B}\Ker \delta_{i-1}$ for all $i\geq 1$. Then we get
the $(A,\,B)$-projective resolution of $_{A}X$ of length $n$

$$0\lra \Img \delta_{n}\lraf{} C_{n-1}\lra\cdots\lra C_1\lra C_0\lra X\lra 0.\quad \quad (\ast_{1})$$

\noindent by Lemma \ref{lemma2.2}. Note that we can write

$$\Img \delta_{n}=\oplus_{j=1}^{m}Q_{j}^{s_{nj}} \mbox{\;and\;} C_{i}=\oplus_{j=1}^{m}Q_{j}^{s_{ij}} ,$$

\noindent where all $s_{ij}$ are non-negative integers for $0\leq i\leq n$ and $1\leq j\leq m$. We claim that $\Img \delta_{n}$ and $C_{i}$ with $0\leq i\leq n-1$ have finite projective dimension. In fact, we may consider the exact sequence
of $A$-modules $0\ra \Ker \delta_{0}\ra C_{0}\lraf{\delta_{0}} X\ra 0$ obtained from $(\ast_{1})$. Since both $X$ and $C_{0}$ have finite projective dimension
by assumption, we have $\pd (\Ker \delta_{0})<\infty$. Then one proceed in the same way from $\Ker \delta_{0}$ in order to show that $C_{1}$ and $\Ker \delta_{1}$ have finite projective dimension, etc. This shows that what we want.
Now by Lemma \ref{lem2.1} we can bound the projective dimension of $_{A}X$:

$$\begin{array}{rl} \pd (_{A}X)&\leq \mbox{max} \{\pd(\Img \delta_{n}),\, \pd(C_{i}),\, i=0,\,\cdots,\,n-1 \}+n\\
&\\
&=\mbox{max} \{\psi(\Img \delta_{n}),\, \psi(C_{i}),\, i=0,\,\cdots,\,n-1 \}+n\\
&\\
&=\mbox{max} \{\psi(\oplus_{j=1}^{m}Q_{j}^{s_{nj}}),\, \psi(\oplus_{j=1}^{m}Q_{j}^{s_{ij}}),\, i=0,\,\cdots,\,n-1 \}+n\\
&\\
&\leq \Psi(\oplus_{j=1}^{m}Q_{j})+n.
\end{array},$$

\medskip
\noindent where $\Psi$ is the {\rm Igusa-Todorov} function. Thus $\fd(A)$ is upper bounded by $\Psi(\oplus_{j=1}^{m}Q_{j})+n$. This completes the proof. $\square$

\medskip
As an immediate consequence of Theorem \ref{thm1}, we have the following.

\begin{Koro} Let $B$ be a subalgebra of an Artin algebra $A$ with the same identity such that $\rad(B)$ is a left ideal in $A$ and $\rad(B)A=\rad(A)$.
If $B$ is representation-finite, then $\fd (A)<\infty$.
\end{Koro}

{\bf Proof.} Consider the inclusion map $i: B\rightarrow A$. Note that, if $\rad(B)$ is a left ideal in $A$ and $\rad(B)A=\rad(A)$, then $\fd(i)\leq \gd(i)\leq1$ by \cite[Proposition 2.19]{XX}. Therefore, if $B$ is representation-finite, by Theorem \ref{thm1}, we have $\fd (A)<\infty$. $\square$

\medskip
As another consequence of Theorem \ref{thm1}, we have the following result.

\begin{Koro} Let $\varphi: B \to A$  be an extension of Artin algebras such that $2\leq\fd(\varphi)<\infty$. Suppose that $\pd(_{B}A)<\infty$ and $A_{B}$ is projective. If $B$ is representation-finite, then $\fd (A)<\infty$. \end{Koro}

{\bf Proof.} By Theorem \ref{thm1}, it suffices to prove that, for any $A$-module $X$ with finite projective dimension, $_{A}A\otimes_{B}X$ has finite projective dimension. Let $X$ be an $A$-module with $\pd(_{A}X)<\infty$. Then, viewing $X$ as a $B$-module, we have $\pd(_{B}X)\leq \pd(_{A}X)+\pd(_{B}A)$. By assumption, we get $\pd(_{B}X)<\infty$, and say $\pd(_{B}X)=m<\infty$. Take a $B$-projective resolution of $_{B}X$ of length $m$
$$0\lra P_{m}\lra P_{m-1}\lra \cdots \lra P_{0}\lra _{B}X\lra 0.$$

\noindent Since $A_{B}$ is projective, the sequence
$$0\lra A\otimes_{B}P_{m}\lra A\otimes_{B}P_{m-1}\lra \cdots \lra A\otimes_{B}P_{0}\lra A\otimes_{B}X\lra 0$$

\noindent is exact and hence an $A$-projective resolution of $_{A}A\otimes_{B}X$, which means that $\pd(_{A}A\otimes_{B}X)<\infty$. $\square$

\bigskip
Remark that, more generally, the above corollary still holds whenever $\fd (B)<\infty$. In fact, let $X$ be an $A$-module with $\pd(_{A}X)<\infty$.
Since $2\leq\fd(\varphi)=n<\infty$, by the proof of Theorem \ref{thm1}, $_{A}X$ has the $(A,\,B)$-projective resolution of length $n$:

$$0\lra Y\lraf{} C_{n-1}\lra\cdots\lra C_1\lra C_0\lra X\lra 0,$$

\noindent such that $_{A}Y$ and $_{A}C_{i}$ with $0\leq i\leq n-1$ have finite projective dimension and that every $C_{i}$ for $0\leq i\leq n-1$ can be expressed the following form:

$$C_{i}=A\otimes_{B}M_{i},$$

\noindent where $_{B}M_{i}$ has finite projective dimension. Assume that $\fd (B)=m<\infty$. Then $_{B}M_{i}$ has a projective resolution of length $m$

$$0\lra Q_{m}\lra Q_{m-1}\lra \cdots \lra Q_{0}\lra _{B}M_{i}\lra 0.$$

\noindent Since $A_{B}$ is projective, the sequence
$$0\lra A\otimes_{B}Q_{m}\lra A\otimes_{B}Q_{m-1}\lra \cdots \lra A\otimes_{B}Q_{0}\lra A\otimes_{B}M_{i}\lra 0$$

\noindent is exact and hence an $A$-projective resolution of $_{A}A\otimes_{B}M_{i}$, which means that $\pd(_{A}A\otimes_{B}M_{i})\leq m$.
Note that $\pd(_{B}Y)\leq \pd(_{A}Y)+\pd(_{B}A)<\infty$, so, by the same way, $\pd(_{A}Y)\leq\pd(_{A}A\otimes_{B}Y)\leq m$. Now we can estimate the projective dimension of $_{A}X$:

$$\begin{array}{rl} \pd (_{A}X)&\leq \mbox{max} \{\pd(Y),\, \pd(C_{i}),\, i=0,\,\cdots,\,n-1 \}+n\\
&\\
&\leq m+n.
\end{array}$$

\noindent This implies that $\fd (A)<\infty$.

\section{Quotient algebras and finitistic dimensions \label{qfd}}

In this section, we shall use representation-theoretical properties of quotient algebras to approach the finiteness of the finitistic dimension of given algebras and prove Theorem \ref{thm2}.

\medskip
{\bf Proof of Theorem \ref{thm2}.} Let $X$ be an $A$-module with finite projective dimension. Consider the exact sequence of $A$-modules
$$0\lra J\Omega_{A}(X)\lra \Omega_{A}(X)\lra\Omega_{A}(X)/J\Omega_{A}(X)\lra 0. $$

\noindent Since $IJ\Omega_{A}(X)\subseteq IJ\rad(P(X))=IJ\rad(A)P(X)\subseteq IJKP(X)=0$ by
assumption, where $P(X)$ is the projective cover of $X$, we have $Y:=J\Omega_{A}(X)$ is an $A/I$-module. Clearly, $Z:=\Omega_{A}(X)/J\Omega_{A}(X)$
is an $A/J$-module.

If $A/I$ and $A/J$ are $A$-syzygy-finite, then there is a non-negative integer $n$, an $A$-module $M$ and an $A$-module $N$ such that $\Omega_{A}^{n}(Y)\in \add (_{A}M)$ and $\Omega_{A}^{n}(Z)\in \add (_{A}N)$. Using Horseshoe Lemma to the above exact sequence, we obtain the
following exact sequence
$$0\lra \Omega_{A}^{n}(Y)\lra \Omega_{A}^{n+1}(X)\oplus P\lra\Omega_{A}^{n}(Z)\lra 0$$

\noindent with $P$ projective $A$-module. Now we can bound the projective dimension of $_{A}X$:

$$\begin{array}{rl} \pd (_{A}X)&\leq \pd (\Omega_{A}^{n+1}(X))+n+1\\
&\\
&=\pd (\Omega_{A}^{n+1}(X)\oplus P)+n+1\\
&\\
&\leq\Psi(\Omega_{A}^{n+1}(Y)\oplus \Omega_{A}^{n+2}(Z))+n+3\\
&\\
&\leq\Psi(\Omega_{A}^{}(M)\oplus \Omega_{A}^{2}(N))+n+3.\\
\end{array},$$

\medskip
\noindent where $\Psi$ is the {\rm Igusa-Todorov} function. Thus $\fd(A)$ is upper bounded by $\Psi(\Omega_{A}^{}(M)\oplus \Omega_{A}^{2}(N))+n+3$ and
$\fd(A)<\infty$.  $\square$

\medskip
The proof of Theorem \ref{thm2} is similar to that of  \cite[Theorem 3.2]{X1} and \cite[Theorem 3.1]{Wjq2}, in which the {\rm Igusa-Todorov} function
is used. However, the difference is that the syzygy shifted sequences is employed in theorem above. It is worth noting that
our result unifies many of the results in literature in this direction, that is, many known results can be obtained from Theorem
\ref{thm2}. In what follows, we shall illustrate it.

If we take $K=A$, then we reobtain \cite[Theorem 3.2]{X1}.

\begin{Koro} (\cite{X1}) Let $A$ be an Artin algebra and $I$, $J$ be two ideals of $A$ such that $IJ=0$. If $A/I$ and $A/J$ are $A$-syzygy-finite (for example, $A/I$ and $A/J$ are representation-finite), then $\fd(A)<\infty$. In particular, algebras with radical-square-zero have finite finitistic dimension.
\end{Koro}

\medskip
If we take $K=\rad (A)$, we have the following result, which recovers \cite[Theorem 3.1]{Wjq2} by \cite[Corollary 2.8]{Wjq2}.

\begin{Koro} Let $A$ be an Artin algebra and let $I$, $J$ be two ideals of $A$ such that $IJ\rad (A)=0$, and that both $_{A}I$ and $_{A}J$ have finite projective dimension. If $A/I$ is $A/I$-syzygy-finite and $A/J$ is $A/J$-syzygy-finite, then the finitistic dimension of $A$ is finite.
\end{Koro}

\medskip
If we set $I=J=\rad^{n}(A)$ and $K=\rad (A)$, we reobtain the main result in \cite{W}.

\begin{Koro} Let $A$ be an Artin algebra with $\rad^{2n+1}(A)=0$. If $A/\rad^{n}(A)$ is $A$-syzygy-finite (for example, $A/\rad^{n}(A)$ is representation-finite), then $\fd(A)<\infty$. In particular, algebras with radical-cube-zero have finite finitistic dimension.
\end{Koro}

\medskip
As other consequence of Theorem \ref{thm2}, we have the following result, generalizing the results of \cite[Lemma 3.6 and Corollary 3.8]{X1} and
\cite[Corollary 3.2 and Proposition 3.5]{Wjq2}.

\begin{Koro} Let $A$ be an Artin algebra with an ideal $I$ such that $A/I$ is $A$-syzygy-finite (for example, $A/I$ is representation-finite).
Then $\fd(A)<\infty$ if one of the following conditions is satisfied: $(1)$\,$I\rad^{2}(A)=0$; $(2)$\,$\rad (A)I\rad (A)=0$;  $(3)$\,$I^{2}\rad (A)=0$.
\end{Koro}

\section{Left idealized extensions and finitistic dimensions \label{lfd}}
In this section, we shall employ left idealized extensions to study the finitistic dimension conjecture and give a proof of Proposition \ref{propo1}.
More precisely, we consider the following question: given a chain of subalgebras of an Artin algebra $A$

$$B=A_{0}\subseteq A_{1}\subseteq\cdots  \subseteq A_{s-1} \subseteq A_{s}=A$$

\noindent such that $\rad(A_{i-1})$ is a left ideal of $A_{i}$ for all $1\leq i\leq s$, if $A$ is representation-finite,
is the finitistic dimension of $B$ finite?

It is known that an affirmative answer to this question will imply that the finitistic
dimension conjecture over the field is true. It was proved by Xi in \cite{X1} that, given such a chain with $s\leq 2$, if $A$ is representation-finite,
then $\fd (B)<\infty$. A natural question is: Is it possible to show that the finitistic dimension of $B$ is finite if $s>2$? In this direction, Wei gave
in \cite[Theorem 2.9]{Wjq} an affirmative answer under some homological conditions. In this section, we shall give a partial answer for the case $s> 2$ by imposing the condition concerning syzygy-finite algebras, which generalizes some results in \cite{X1,WX}.

Let us start with the following two lemmas from \cite{WX,X2}, which establish a way of lifting modules over a subalgebra to modules over its extension algebra.

\begin{Lem} (\cite[Lemma 3.2]{X2}) \label{lem5.1} Let $A$ be an Artin algebra and $B$ be a subalgebra of $A$ with the same identity such that $\rad (B)$ is a left ideal of $A$. Then, for any $B$-module $X$, $\Omega_{B}^{i}(X)$ is a torsionless $A$-module for all $i\geq 2$ and there is a projective $A$-module $P$ and an $A$-module $Y$ such that
$\Omega_{B}^{i}(X)\simeq \Omega_{A}(Y)\oplus P$ as $A$-modules. \end{Lem}

\begin{Lem} (\cite[Lemma 3.5]{WX}) \label{lem5.2} Let $A$ be an Artin algebra and $B$ be a subalgebra of $A$ with the same identity. Suppose that $I$ is an ideal of $B$ and is also a left ideal of $A$. Then, for any torsionless $B$-module $X$, $IX$ is a torsionless $A$-module and there is a projective $A$-module $Q$ and an $A$-module $Z$ such that
$IX\simeq \Omega_{A}(Z)\oplus Q$ as $A$-modules. \end{Lem}

\medskip
In the following, we employ the finitistic dimension of bigger algebras to approach that of smaller algebras for a left idealized extension.

\begin{Lem} Let
$$B=A_{0}\subseteq A_{1}\subseteq\cdots  \subseteq A_{s-1}\subseteq A_{s}=A$$

\noindent be a chain of subalgebras of an Artin algebra $A$ such that $I_{i-1}$ is an ideal of $A_{i-1}$ and is also a left ideal of $A_{i}$ for all $1\leq i\leq s$ with $s$ being a positive integer. If $A$ is 1-syzygy-finite and $B/I_{s-1}\cdots I_{1}I_{0}$ is $B$-syzygy-finite (for example, $B/I_{s-1}\cdots I_{1}I_{0}$ is representation-finite), then $\fd(B)<\infty$. \label{lem5.3}\end{Lem}

{\bf Proof}. First observe that $I:=I_{s-1}\cdots I_{1}I_{0}$ is an ideal of $B$ contained in $I_{0}$ by assumption. Let $X$ be a $B$-module with finite projective dimension. Then we can form an exact sequence of $B$-modules:

$$0\lra I\Omega_{B}(X)\lra \Omega_{B}(X)\lra\Omega_{B}(X)/I\Omega_{B}(X)\lra 0.$$

\noindent Clearly, $Y:=\Omega_{B}(X)/I\Omega_{B}(X)$ is a $B/I$-module and hence there is a non-negative integer $n$ and a $B$-module $M$ such that
$\Omega_{B}^{n}(Y)\in \add (_{B}M)$, since $B/I$ is $B$-syzygy-finite. Note that $\Omega_{B}(X)$ is a torsionless $B$-module, so $I_{0}\Omega_{B}(X)$ is a torsionless $A_{1}$-module by Lemma \ref{lem5.2}. Inductively, by Lemma \ref{lem5.2} again, we obtain that $Z:=I\Omega_{B}(X)$ is a torsionless $A$-module. Hence, there is a projective $A$-module $Q$ and an $A$-module $W$ such that $Z\simeq \Omega_{A}(W)\oplus Q$ as $A$-modules. Since $A$ is 1-syzygy-finite, there exists an $A$-module $N$ such that $\Omega_{A}^{}(W)\in \add (_{A}N)$, which means that $Z\in\add (_{A}N\oplus A)$.

Taking the $n$-th syzygy of the above exact sequence, by Horseshoe Lemma, we obtain an exact sequence of $B$-modules

$$0\lra \Omega_{B}^{n}(Z)\lra \Omega_{B}^{n+1}(X)\oplus P\lra\Omega_{B}^{n}(Y)\lra 0$$

\noindent with $P$ projective $B$-module. Now we can bound the projective dimension of $_{B}X$:

$$\begin{array}{rl} \pd (_{B}X)&\leq \pd (\Omega_{B}^{n+1}(X))+n+1\\
&\\
&=\pd (\Omega_{B}^{n+1}(X)\oplus P)+n+1\\
&\\
&\leq\Psi(\Omega_{B}^{n+1}(Z)\oplus \Omega_{B}^{n+2}(Y))+n+3\\
&\\
&\leq\Psi(\Omega_{B}^{n+1}(N\oplus A)\oplus \Omega_{B}^{2}(M))+n+3\\
\end{array},$$

\medskip
\noindent where $\Psi$ is the {\rm Igusa-Todorov} function. Thus $\fd(B)$ is upper bounded by $\Psi(\Omega_{B}^{n+1}(N\oplus A)\oplus \Omega_{B}^{2}(M))+n+3$. This completes the proof. $\square$

\medskip
Note that Lemma \ref{lem5.3} recovers \cite[Theorem 3.1]{X1} if we take $s=1$. The next result is a variation of Lemma \ref{lem5.3}.

\begin{Lem} Let

$$B=A_{0}\subseteq A_{1}\subseteq\cdots  \subseteq A_{s-1}\subseteq A_{s}=A$$

\noindent be a chain of subalgebras of an Artin algebra $A$ such that $I_{i-1}$ is an ideal of $A_{i-1}$ and is also a left ideal of $A_{i}$ for all $1\leq i\leq s$ with $s$ being a positive integer. Suppose that $I_{0}$ is the Jacobson radical $\rad (B)$ of $B$. If $A$ is 1-syzygy-finite and $A_{1}/I_{s-1}\cdots I_{1}$ is $B$-syzygy-finite (for example, $A_{1}/I_{s-1}\cdots I_{1}$ is representation-finite), then $\fd(B)<\infty$. \label{lem5.4}\end{Lem}

{\bf Proof}. Given a $B$-module $X$ with finite projective dimension, we consider $\Omega_{B}^{2}(X)$ instead of $\Omega_{B}^{}(X)$.
By Lemma \ref{lem5.1}, $\Omega_{B}^{2}(X)$ is a torsionless $A_{1}$-module. Then we can form an exact sequence of $B$-modules:
$$0\lra I_{s-1}\cdots I_{1}\Omega_{B}^{2}(X)\lra \Omega_{B}^{2}(X)\lra\Omega_{B}^{2}(X)/I_{s-1}\cdots I_{1}\Omega_{B}^{2}(X)\lra 0.$$

\noindent By the argument in the proof of Lemma \ref{lem5.3} we obtain the lemma. $\square$

\medskip
Here, we understand $A_{1}/I_{s-1}\cdots I_{1}=0$ if $s=1$, which means that $A_{1}/I_{s-1}\cdots I_{1}$ being $B$-syzygy-finite always holds. Let us remark that Lemma \ref{lem5.4} recovers \cite[Theorem 3.1]{X1} if we take $s=1$, and extends \cite[Corollary 3.10]{WX} if we take $s=2$.

Combining Lemma \ref{lem5.3} with Lemma \ref{lem5.4}, we prove Proposition \ref{propo1}. As an immediate consequence of Proposition \ref{propo1}, we have the following corollary, which is a partial answer to the question by Xi in his website~(see http://math0.bnu.edu.cn/~ccxi/Problems.php).

\begin{Koro} Let $D\subseteq C\subseteq B\subseteq A$
be a chain of subalgebras of an Artin algebra $A$ such that the radicals of $D$, $C$ and $B$ are left ideals of $C$, $B$
and $A$, respectively. Suppose that $A$ is 1-syzygy-finite. If either $D/\rad (B)\rad (C)\rad (D)$ or $C/\rad (B)\rad (C)$ is $D$-syzygy-finite,
then $\fd(D)<\infty$. In particular, if either $D/\rad (B)\rad(C)\rad (D)$ or $C/\rad (B)\rad (C)$ is representation-finite, then $\fd(D)<\infty$.
\end{Koro}

We end this section with an example showing that our results do apply to check the finiteness of the finitistic dimension of some algebras.

\bigskip
\noindent{\bf Example 1.}\,\,(\cite{WX})\,\,Let $A$ be the algebra given by the quiver with relation:

$$\xymatrix{\circ
&\circ\ar[l]_{\lambda}^(1){5}^(0){2}&\circ\ar[l]_{\epsilon}^(0){3}&\circ\ar[l]_{\xi}^(0){1} &\circ\ar[l]_{\beta}^(){}^(0){4}&\circ\ar[l]_{\alpha}^(){}^(0){6}},\,\,\,\,\alpha\beta\xi\epsilon\lambda=0.$$

\noindent Then $A$ is a Nakayama algebra and hence 1-syzygy-finite. Let $C$ be the subalgebra of A generated by the set $\{e_{1}, e_{2'}:=e_{2}+e_{4}+e_{5}, e_{3'}:=e_{3}+e_{6}, \lambda, \beta, \alpha+\epsilon, \gamma:=\xi\epsilon, \delta:=\beta\xi\}$, which is given by the quiver with relations:

$$
\xymatrix{\circ\ar@<0.4ex>[r]^{\gamma}
&\circ\ar@<0.4ex>[l]^{\beta}^(1){1}^(0){2'}\ar@(ur,ul)_{\lambda}\ar@<0.4ex>[r]^{\delta}
&\circ\ar@<0.4ex>[l]^{\alpha+\epsilon}^(0){3'}},\quad \beta\gamma=\delta(\alpha+\epsilon), \gamma\beta=\gamma\delta=\lambda^{2}=\lambda\beta=\lambda\delta=(\alpha+\epsilon)\beta\gamma\lambda=0.
$$

\noindent It is not hard to see that $\rad^{3}(C)\neq 0$ and $\ell\ell^{\infty}(C)=4$, where $\ell\ell^{\infty}(C)$ denotes the infinite-layer length of $C$ (\cite{HLM}). Note also that $C$ is neither a monomial algebra nor a special biserial algebra.
It was proved in \cite{WX} that the finitistic dimension of $C$ is finite. Here, we shall use left idealized extensions to reobtain the finiteness of the finitistic dimension of $C$, though $\rad(C)$ is not a left ideal of $A$. In fact, let $B$ be the subalgebra of A generated by the set $\{e_{1}, e_{2'}:=e_{2}+e_{4}+e_{5}, e_{3'}:=e_{3}+e_{6}, \lambda, \beta, \alpha, \epsilon, \gamma:=\xi\epsilon, \delta:=\beta\xi\}$. Then $B$ is given by the following quiver with relations:

$$
\xymatrix{\circ\ar@<0.4ex>[r]^{\gamma}
&\circ\ar@<0.4ex>[l]^{\beta}^(1){1}^(0){2'}\ar@(ur,ul)_{\lambda}\ar@<1.3ex>[r]^{\delta}
&\circ\ar[l]|{\epsilon}^(0){3'}\ar@<1.3ex>[l]^{\alpha}},\quad \beta\gamma=\delta\epsilon, \gamma\beta=\gamma\delta=\lambda^{2}=\lambda\beta=\lambda\delta=\delta\alpha=\epsilon\beta=\epsilon\delta=\alpha\lambda=\alpha\beta\gamma\lambda=0.
$$

\medskip
\noindent It is easy to check that $C\subseteq B\subseteq A$
is a chain of subalgebras of $A$ such that $\rad(C)$ is a left ideal of $B$ and $\rad(B)$ is a left ideal of $A$. So we have $\fd(C)<\infty$ by
Proposition \ref{propo1}.

\bigskip
%{\bf Acknowledgements.} The author is greatly indebted to Prof. Changchang Xi for his guidance, and thanks the referee for helpful comments and suggestions.
% The research work is partially supported by the NSFC (No. 11331006) and Fundamental Ability Enhancement Project for Young and Middle-aged University Teachers in Guangxi Province (No. 2017KY0860).

\end{document}